\documentclass[11pt]{amsart}
\usepackage{latexsym,amssymb,amsmath,amscd,amsthm}
\numberwithin{equation}{section}
\theoremstyle{remark}
\newtheorem{theorem}{{\bf THEOREM}}[section]

\newtheorem{corollary}{{\bf COROLLARY}}[section]
\newtheorem{example}{{\bf EXAMPLE}}[section]
\newtheorem{proposition}{{\bf PROPOSITION}}[section]
\newcommand{\bq}{\begin{equation}}
\newcommand{\bea}{\begin{array}}
\newcommand{\eea}{\end{array}}

\newcommand{\ga}{\alpha}

\newcommand{\gD}{\Delta}
\newcommand{\gl}{\lambda}

\newcommand{\gb}{\beta}

\newcommand{\ul}[1]{\underline{#1}}

\newcommand{\gG}{\Gamma}

\newcommand{\pp}{\partial}

\newcommand{\tl}{\tilde}

\newcommand{\nm}{\left[\begin{array}{c}
n\\
m\end{array}\right]}
\newcommand {\nk}{\left(\begin{array}{c}
n\\
k\end{array}\right)}

\newcommand{\bs}{\blacksquare}

\newcommand{\bgs}{\bigstar}

\newcommand{{\DDD}}{D\!\!\!\!\!\!-}
\title{REMARKS ON Q-CALCULUS AND INTEGRABILITY}
\author{Robert Carroll\\University of Illinois, Urbana, IL 61801}
\date{August, 2002\thanks{email: rcarroll@math.uiuc.edu}}
\begin{document}

\bibliographystyle{plain}

\maketitle

%\newpage
%\tableofcontents
%\addcontentsline{toc}{section}{subsection}
%\setcounter{tocdepth}{3}

\section{INTRODUCTION}
\renewcommand{\theequation}{1.\arabic{equation}}
\setcounter{equation}{0}

Integrable systems such as qKP have been frequently studied in recent years
(see e.g. \cite{a2,a3,cxx,c1,c2,f1,f2,h1,i6,k3,k2,t1,t2}).  Also various
noncommutative (NC) integrable models connecting frequently to Moyal
deformations arise in the literature (see e.g. \cite{ch,d1,d2,d3,j1,j2}).  In
view of the intimate relations between integrable systems, differential
calculi, and virtually all aspects of theoretical physics (cf.
\cite{ch,cxx,c4}), as well as the profound connections between quantum groups
(QG) and quantum field theory (QFT) for example (cf. \cite{b1,b2}), it seems
compelling to pursue further the relations between QG, integrable systems,
differential calculi, and zero curvature (cf. \cite{cxx,c3,d1,d2,d3,d4,j1,j2}).
The present article is primarily concerned with various forms of Hirota
equations.

\section{SOME BACKGROUND}
\renewcommand{\theequation}{2.\arabic{equation}}
\setcounter{equation}{0}

We recall first the derivation of KP from a differential calculus in 
\cite{cxx,c3,d3,d4} (fundamentals of QG theory will be assumed as we proceed
(cf. \cite{cxx,k1,m1}).
\begin{example}
Consider a calculus based on {\bf (A)}  $dt^2=dx^2=dxdt+dtdx=0$  {\bf (B)}
$[dt,t]=[dx,t]=[dt,x]=0$ and $[dx,x]=\eta dt$.  Assumingt the Leibnitz rule
$d(fg)=(df)g+f(dg)$ for functions and $d^2=0$ one obtains ${\bf (A1)}\,\,
df=f_xdx+(f_t+(1/2)\eta f_{xx})dt$.  For a connection $A=wdt+udx$ the zero
curvature condition $F=dA+A^2=0$ leads to ${\bf
(A2)}\,\,(u_t-w_x+(\eta/2)u_{xx}+\eta uu_x=0$ which for $w_x=0$ is a form of
Burger's equation.
$\hfill\bs$
\end{example}
\begin{example}
Next consider {\bf (A)}  $[dt,t]=[dx,t]=[dt,x]=[dy,t]=[dt,y]=[dy,y]=0$ with
{\bf (B)}  $[dx,x]=2bdy$ and $[dx,y]=[dy,x]=3adt$.  Further {\bf (C)}  
$dt^2=dy^2=dtdx+dxdt=dydt+dtdy=dydx+dxdy=0$.  Then ${\bf (A3)}\,\, df=f_xdx+
(f_y+bf_{xx})dy+(f_t+3af_{xy}+abf_{xxx}dt)$.  For $A= vdx+wdt+udy$ one finds
that $dA+A^2=F=0$ implies
\bq\label{1}
u_x=v_y+bv_{xx}+2bvv_x;\,\,w_x=3av_{xy}+abv_{xxx}+3auv_x+3av(v_y+bv_{xx});
\end{equation}
$$w_x+bw_{xx}=u_t+3au_{xy}+abu_{xxx}+3auu_x-v[2bw_x-3a(u_y+bu_{xx})]$$
Taking e.g. $w_x=(3a/2b)u_y+(3a/2)u_{xx}$ in the last equation to decouple one
arrives at ${\bf (A4)}\,\,\pp_x(u_t-(ab/2)u_{xxx}+3auu_x)=(3a/2b)u_{yy}$;
for suitable $a,b$ this is KP.
$\hfill\bs$
\end{example}
\indent
It is surprisingly difficult to convert these examples into meaningful 
q-calculus equations and in that spirit for guidance we were motivated to
develop many formulas concerning qKP, qKdV, etc.  It is convenient and
hopefully useful to write down also various formulas involving q-calculus. 
Thus the presentation has a sort of guide book spirit in which 
a number of formulas are
extracted from various references and some new formulas are derived
(especially in the context of Hirota type formulas).  We recall first from
\cite{a2,a3,cxx,f1,f2,h1,i6,k3,k2,t1} that for ${\bf (A5)}\,\,D_qf(x)=
[f(qx)-f(x)]/(q-1)x$ and $Df(x)=f(qx)$ one writes e.g. ${\bf (A6)}\,\,
L=D_q+a_0+\sum_1^{\infty}a_iD_q^{-i}$ and qKP is defined via ${\bf (A7)}\,\,
\pp_jL=[L_{+}^j,L]$.  Note here also $L=SD_qS^{-1}$ where ${\bf (A8)}\,\,
S=1+\sum_1^{\infty}\tl{w}_jD_q^{-j}$ with $\pp_jS=-L_{-}^jS$ (recall one often
takes $t=(t_1,t_2,\cdots)$ with $\pp_j\sim \pp/\pp t_j$ and it is
often convenient to insert x separately instead of identifying $t_1=x$).  Now
e.g. write
\bq\label{2}
exp_q(x)=\sum_0^{\infty}\frac{(1-q)^kx^k}{(q;q)_k}=exp\left(\sum_1^{\infty}
\frac{(1-q)^kx^k}{k(1-q^k)}\right)
\end{equation}
where $(q;q)_k=(1-q)\cdots(1-q^k)$ (nonstandard definition).
Then $D_qexp_q(xz)=zexp_q(xz)$ and one
defines q-Schur polynomials $\tl{p}_k$ via
\bq\label{3}
\sum_{{\bf Z}}\tl{p}_k(x,t_1,\cdots)z^k=exp_q(xz)exp(\sum_1^{\infty}t_kz^k)
\end{equation}
Then $D_q\tl{p}_k=\tl{p}_{k-1}=\pp_1\tl{p}_k\,\,(\pp_1=\pp/\pp t_1)$ and
$\tl{p}_k=0\,\,(k<0),\,\,
\tl{p}_0=1,\,\,\tl{p}_1=x+t_1,\,\,\tl{p}_2=[(q-1)/(q^2-1)]x^2+t_1x+
t_1^2/2+t_2,\cdots$. 
The  qKP wave functions are defined via ($\xi=\sum_1^{\infty}t_kz^k=\xi(t,z)$)
\bq\label{4}
\psi_q=S
exp_q(xz)exp(\xi);\,\,\psi_q^*=(S^*)^{-1}_{\frac{x}{q}}exp_{\frac{1}{q}}(-xz)exp
(-\xi)
\end{equation}
for $S=1+\sum_1^{\infty}\tl{w}_iD_q^{-i}$,
where for $V=\sum_1^{\infty}v_iD_q^{-i}$, one defines $V_{\frac{x}{q}}=\sum
v_i(x/q)q^iD_q^{-1}$ (note $e_q(xz)^{-1}=e_{1/q}(-xz)$).
It follows that $L\psi_q=z\psi_q$ and $L^*_{x/q}\psi_q^*
=z\psi^*_q$ with a q-Hirota bilinear identity ${\bf (A9)}\,\,Res D_q^n\pp^{\ga}
\psi_q\psi_q^*=0$ where $\pp^{\ga}\sim\pp_1^{\ga_1}\pp_2^{\ga_2}\cdots$.
This is proved by an extended lemma of Dickey (cf. \cite{i6}) but we will
prefer to use somewhat different forms below as well as some new variations.
The q-tau functions are determined via ${\bf (A10)}\,\,\tau_q(x,t)=
\tau(t+c(x))$ where 
$$c(x)=[x]_q=\left(\frac{(1-q)x}{(1-q)},\frac{(1-q)^2x^2}{2(1-q^2)},
\frac{(1-q)^3x^3}{3(1-q^3)},\cdots\right)$$
There results
\bq\label{5}
\psi_q=\frac{\tau_q(x,t-[z^{-1}])}{\tau_q(x,t)}exp_q(xz)exp(\xi);\,\,\psi_q^*=
\frac{\tau_q(x,t+[z^{-1}])}{\tau_q(x,t)}exp_{\frac{1}{q}}(-xz)exp(-\xi)
\end{equation}
where $[z]=(z,z^2/2,z^3/3,\cdots)$.  Vertex operators for KP are defined via
$$X(t,z)=exp(\xi)exp(-\sum_1^{\infty}(\pp_i/i)z^{-i})=exp(\xi(t,z)exp(-\xi(
\tl{\pp},z^{-1}))$$
where $\tl{\pp}\sim (\pp_1,\pp_2/2,\pp_3/3,\cdots)$.  Then
for KP
\bq\label{6}
X_q=e_q(xz)X(t,z);\,\,\tl{X}_q=e_q(xz)^{-1}X(-t,z)=e_{1/q}(-xz)exp(-\xi(t,z)
exp(\tl{\pp},z^{-1})
\end{equation}
We recall also another version of {\bf (A9)} in the form ($n>m$)
\bq\label{7}
\oint_{\infty}D^n(X_q(x,t,z)\tau_q(t))D^{m+1}(\tl{X}_q(x,t',z)\tau_q(t'))dz=0
\end{equation}
\indent
{\bf REMARK 2.1.}  For qKdV one has ${\bf (A11)}\,\,L^2=D+q^2+(q-1)xuD_q+u$
and $u=D_q\pp_1log(\tau_q(x,t)\tau_q(qx,t))$.  We will say more about this
in Section 5.$\hfill\bs$
\\[2mm]\indent
{\bf REMARK 2.2.}  The standard Schur polynomials are defined via
\bq\label{8}
p_n(y)=\sum\left(\frac{y_1^{k_1}}{k_1!}\right)\left(\frac{y_2^{k_2}}
{k_2!}\right)\cdots;\,\,\,\sum_0^{\infty}p_n(y)z^n=exp(\xi(y,z))
\end{equation}
Further with the prescription $u=\pp^2_1log \tau$ of standard KP theory one
has from the bilinear identity (with $x=t_1$ here and $\pp\sim\pp_1$)
\bq\label{9}
\pp_{n-1}\pp^{-1}u=\frac{1}{\tau^2}p_n(\tl{\pp})\tau\cdot\tau=s_n;\,\,
K_4=\pp s_4=\frac{1}{4}\pp^3u+3u\pp u+\frac{3}{4}\pp^{-1}\pp_2^2u=\pp_3u
\end{equation}
(see \eqref{10} for the bilinear notation $\tau\cdot\tau$) .  Here the Hirota
equations (see below for a derivation) are
\bq\label{10}
\left(\sum_0^{\infty}p_n(-2y)p_{n+1}(\tl{\pp})exp(\sum_1^{\infty}y_i\pp_i)
\right)\tau\cdot\tau=0;\,\,
\pp_j^ma\cdot b=(\frac{\pp^m}{\pp
s_j^m})a(t_j+s_j)b(t_j-s_j)|_{s=0}
\end{equation}
The corresponding formula for qKP will be written down later.  Note the
bilinear form of KP is ${\bf
(A12)}\,\,(\pp^4+3\pp^2_2-4\pp_1\pp_3)\tau\cdot\tau=0$.
$\hfill\bs$
\\[2mm]\indent
{\bf REMARK 2.3.}  We collect a few more formulas.  First note ${\bf (A13)}\,\,
D_q^*=-(1/q)D_{1/q}=-D_qD^{-1}$ since
$D_qD^{-1}f(x)=D_qf(q^{-1}x)=[f(x)-f(q^{-1}x)]/(q-1)x=q^{-1}D_{1/q}f$.  One
also notes for $\xi=\sum t_iz^i$ that ${\bf
(A14)}\,\,\psi_q=S(exp_q(xz)exp(\xi))=(1+\sum\tl{w}_iz^{-i})exp_q(xz)exp(\xi)=
\hat{\psi}_qexp_q(xz)exp(\xi)$ (note $\tl{w}_i$ depends on q - cf. {\bf (A30)}
below).  Further ${\bf (A15)}\,\,V=\sum v_iD_q^i\Rightarrow
V^*=\sum(D^*_q)^iv_i$ where $D^*_q=-q^{-1}D_{1/q}$ and
$(S^*)^{-1}$ is written as $(S^*)^{-1}=1+\sum\tl{w}^*_iD_q^{-i}$. 
Note also
${\bf
(A16)}\,\,L\psi_q=z\psi_q,\,\,\pp_n\psi_q=L^n_{+}\psi_q,\,\,L^*_{x/q}\psi^*_q
=z\psi^*_q$ and $\pp_n\psi_q^*=-(L^*_{x/q})^n_{+}\psi_q^*$.  Further
$$\psi_q^*=(S^*)_{x/q}^{-1}exp_{1/q}(-xz)exp(-\xi)=(1+\sum\tl{w}_i^*z^{-i})
exp_{1/q}(-xz)exp(-\xi)$$  
Additional formulas are $D_qD^k(f)=
q^kD^k(D_qf),\,\,D_qf=(D_qf)+D(f)D_q,\,\,D_qD_q^{-1}=D_q^{-1}D_qf=f,$ and thus
${\bf (A17)}\,\,D_q^{-1}f=\sum_{k\ge 0}(-1)^kq^{-k(k+1)/2}D^{-k-1}(D_q^kf)D_q^
{-k-1}$.  This leads to
\bq\label{11}
D_q^nf=\sum_{k\geq 0}\nm_qD^{n-m}(D_q^kf)D_q^{n-m}\,\,\,(n\in{\bf Z})
\end{equation}
where $\nm_q=[n]_q\cdots [n-m+1]_q/[1]_q\cdots [m]_q$ with $[n]_q=(q^n-1)/
((q-1)$
and one notes also
\bq\label{122}
D_q^n=\frac{q^{-n(n-1)/2}}{x^n(q-1)^n}\sum_0^n(-1)^mq^{m(m-1)/2}\nm_q D^{n-m}
\end{equation}
More formulas will appear as we go along.$\hfill\bs$

\section{CLASSICAL HIROTA FORMULAS}
\renewcommand{\theequation}{3.\arabic{equation}}
\setcounter{equation}{0}

We would like to have e.g. a formula analogous to \eqref{10} for qKP, as well
as other ``generic" formulas in the q-context.  In this direction we recall
some results from \cite{a4,a5,ch,t3,t1}.  First note
${\bf
(A18)}\,\,\sum_0^{\infty}p_k(-\tl{\pp})z^{-k}=exp(\sum_1^{\infty}
(-\pp_m/m)z^{-m}$ (recall $\tl{\pp}=(\pp_1,\pp_2,\cdots)$ and one writes
\bq\label{15}
\psi(t,\mu)=X(t,\mu)\tau/\tau;\,\,X(t,\mu)=e^{\xi(t,\mu)}exp(-\sum_1^{\infty}
(\pp_m/m)z^{-m})
\end{equation}
so formally
\bq\label{16}
\psi(t,\mu)=e^{\xi(t,\mu)}\sum_0^{\infty}p_k(-\tl{\pp})\mu^{-k}\tau/\tau=
e^{\xi(t,\mu)}\frac{\tau(t-[\mu^{-1}])}{\tau}
\end{equation}
with ${\bf
(A19)}\,\,\tau(t-[\mu^{-1}])=(\sum_0^{\infty}p_k(-\tl{\pp})\mu^{-k})
\tau$.  Similarly ${\bf (A20)}\,\,\psi^*(t,\gl)=X^*(t,\gl)\tau/\tau$ with
$X^* (t,\gl)=exp(-\xi(t,\gl)exp(\sum_1^{\infty}(\pp_m/m)\gl^{-m})$ leading
to a natural formula ${\bf (A21)}\,\,\psi^*(t,\gl)=exp(-\xi(t,\gl)\sum
p_k(\tl{\pp})\gl^{-k}\tau/\tau= exp(-\xi)\tau(t+[\gl^{-1}])/\tau$.  Note
that in \cite{a5}, instead of the more standard W with $W=1+\sum w_j\pp^{-j}$
satisfying $L=W\pp W^{-1}$ with $\psi=Wexp(\xi)$ (cf. \cite{ch,c4}), one
works with $\ul{t}= (t_1+x,t_2,t_3,\cdots)$ and sets ${\bf
(A22)}\,\,\psi(x,t,z)=(Sexp(xz)exp(\xi))=exp(xz)exp(\xi)\tau(\ul{t}-
[z^{-1}])/\tau(t)$.  Then explicitly it is known that 
\bq\label{18}
S=\sum_0^{\infty}\frac{p_n(-\tl{\pp})\tau(\ul{t})}{\tau(\ul{t})}\pp^{-n}
\end{equation}
(standard Schur polynomials).  In this spirit one has also
\bq\label{19}
\psi^*=(S^*)^{-1}e^{-xz}e^{\xi}=e^{-xz-\xi}\frac{\tau(\tl{t}+[z^{-1}])}
{\tau(\ul{t})}
\end{equation}
(note also in \cite{ch} ${\bf (A23)}\,\,W=1+\sum w_j\pp^{-j},
\,\,W^*=1+\sum w^*_j\pp^{-j}\,\,w_j=(1/\tau)p_j(-\tl{\pp})\tau,$
and $w_j^*=(1/\tau)p_j(\tl{\pp})\tau$). 
Clearly the operator $exp(\sum
(-\pp_i/i)z^{-i})$ acting on $\tau$ simply translates variables and $t_1$ or
$t_1+x$ is equally affected; one chooses $\tau(\ul{t}-[z^{-1}])$ by
``design". Next we recall the classical differential Fay identity
\bq\label{21}
\pp_1\tau(t-[\mu^{-1}])\tau(t-[\gl^{-1}])-\tau(t-[\mu^{-1}])\pp_1
\tau(t-[\gl^{-1}])-
\end{equation}
$$-(\gl-\mu)\tau(t-[\mu^{-1}])\tau(t-[\gl^{-1}])+(\gl-\mu)\tau(t)\tau
(t-[\mu^{-1}]-[\gl^{-1}])=0$$
and write this as in 
\cite{t3}
\bq\label{22}
\frac{\tau(t-[\mu^{-1}]-[\gl^{-1}])\tau(t)}{\tau(t-[\mu^{-1}])
\tau(t-]\gl^{-1}])}=1+\frac{1}{\mu-\gl}\left[\pp_1log(t-[\mu^{-1}])-
\pp_1log(t-[\gl^{-1}])\right]
\end{equation}
This can be used to generate a useful formula in dealing with dispersionless
Hirota equations (cf. \cite{ch,c4,t3} and see remarks below).
In connection with Hirota formulas we recall also for $t\ne t'$
\bq\label{24}
\oint \psi^*(t,\gl)\psi(t',\gl)d\gl=0\leadsto \oint \tau(t+[\gl^{-1}])
\tau(t'-[\gl^{-1}])e^{\sum (t'_i-t_i)\gl^i}d\gl=0
\end{equation}
Let $t\to t+y,\,\,t'\to t-y$ to obtain
\bq\label{25}
0=\oint \tau(t+y+[\gl^{-1}])\tau(t-y-[\gl^{-1}])e^{-2\sum y_i\gl^i}d\gl=
\end{equation}
$$=\oint e^{\sum y_i\pp_i+\sum\tl{\pp}_i\gl^{-i}}\tau\cdot\tau e^{-2\sum
y_i\gl^i}d\gl
=\oint \sum p_n(-2y)\gl^n\sum p_{\ell}(\tl{\pp})\gl^{-\ell}e^{\sum y_i\pp_i}
\tau\cdot\tau d\gl$$
This leads to the well known (cf. \cite{ch} and \eqref{10}) formula involving
the coefficient of $y_n$
\bq\label{226}
\pp_1\pp_n\tau\cdot\tau=2p_{n+1}(\tl{\pp})\tau\cdot\tau
\end{equation}
Note for the coefficient of $y_n$ one takes first the residue term for $p_n$
where $y_n$ appears bare and then the term for $n=0$ and $\ell=1$ to get
$p_1\sim\pp_1$ and $\pp_n$ from the exponent.
\\[3mm]\indent
Now write from \cite{ch} ${\bf (A24)}\,\,\psi\psi^*=\sum_0^{\infty}s_n
\gl^{-n}$ and $s_n=\sum_0^nw_jw^*_{n-j}$ where $w_0=1$ and $s_0=1$.  Given
\eqref{18} and remarks after {\bf (A23)} we can compute $s_n$ directly and
then one can utilize \eqref{9}.  Thus 
\bq\label{26}
s_n=\sum_0^np_j(-\tl{\pp})\tau p_{n-j}(\tl{\pp})\tau/\tau^2
\end{equation}
This is direct and avoids the Hirota bilinear formalism
(evidently one has $s_0=1, s_1=0, s_2= \pp_1^2log\tau$, etc. in
agreement with \eqref{9}). This also yields an apparently
new formula.
\begin{proposition}
\bq\label{27}
\sum_0^np_j(-\tl{\pp})\tau p_{n-j}(\tl{\pp})\tau=p_n(\tl{\pp})\tau\cdot
\tau
\end{equation}
\end{proposition}
\indent
Now to clarify the ``evaluation" of \eqref{22} we note the important property
\bq\label{28}
log\tau(t-[\gl^{-1}])=log exp\left(\sum(-\pp_i/i)\gl^{-i}\right)\tau
=log\sum
p_k(-\tl{\pp})\gl^{-k}\tau=
\end{equation}
$$=exp\left(\sum(-\pp_i/i)\gl^{-i}\right)log\tau
=\sum p_k(-\tl{\pp})\gl^{-k}log\tau$$
Now for the left side
of \eqref{22} take logarithms to obtain
\bq\label{29}
\sum p_k(-\tl{\pp})\gl^{-k}p_m(-\tl{\pp})\mu^{-m}log\tau+log\tau-
\sum p_m(-\tl{\pp})\mu^{-m}log\tau-
\end{equation}
$$- \sum p_k(-\tl{\pp})\gl^{-k}log\tau =\sum_1^{\infty}p_k(-\tl{\pp})p_m
((-\tl{\pp})\gl^{-k}\mu^{-m}log\tau$$
The terms with $k=m=0,\,\,k=0,$ and $m=0$ combine with $log \tau$ to produce 
zero.  This yields then
\bq\label{30}
\sum_1^{\infty} p_k(-\tl{\pp})p_m(-\tl{\pp})\gl^{-k}\mu^{-m}log\tau=
log\left(1+\sum_1^{\infty}\frac{\mu^{-n}-\gl^{-n}}{\mu-\gl}p_n(-\tl{\pp})
\pp_1log\tau\right)
\end{equation}
One can push this a little further in the spirit of \cite{c4}.  Thus, setting
$\pp_1log\tau(t-[\gl^{-1}])=f(t,\gl)$, we have ${\bf (A25)} \,\,[f(t,\mu)-
f(t,\gl)]/(\mu-\gl)\to\pp_{\gl}f(t,\gl)$ as $\mu\to \gl$, so from \eqref{30}
we obtain
\bq\label{31}
\sum_1^{\infty}p_k(-\tl{\pp})p_m(-\tl{\pp})\gl^{-k-m}log\tau=L=log
(1+\pp_{\gl}\pp_1log\tau(t-[\gl^{-1}])=R=
\end{equation}
$$=log(1+\sum_1^{\infty}\pp_nf(t,\gl)\gl^{-n-1}=log(1+\sum_1^{\infty}
\gl^{-n-1}\pp_n\pp_1log\tau(t-[\gl^{-1}]))=log(D)=$$
$$=log(1+\sum_1^{\infty}\gl^{-n-1}\sum_0^{\infty}p_{\ell}(-\tl{\pp})\gl^{-\ell}
\pp_n\pp_1log\tau)=log(1+\sum \gl^{-n-\ell-1}p_{\ell}(-\tl{\pp})\pp_n\pp_1log
\tau)$$
($1\leq n<\infty$ and $0\leq \ell<\infty$).  Note also we could write
$(\pp_1L)D=\pp_1D$ to get rid of the logarithm; then use $\pp_1^2log\tau=u$
and one has a form of Hirota equation generating function in terms of u alone!
However it is rather too complicated for computations (see below however).
We could also try to duplicate further the procedure of \cite{c4}
(done below).
\\[3mm]\indent
To indicate the direct computations based on L,R,D write out
$(\pp_1L)D=\pp_1D$ in the form
\bq\label{32}
(\sum_1^{\infty} p_k(-\tl{\pp})p_m(-\tl{\pp})\gl^{-k-m}\pp_1log\tau)(1+
\sum_1^{\infty}\gl^{-n-1}\sum_0^{\infty}p_{\ell}(-\tl{\pp})\gl^{-\ell}
\pp_n\pp_1log\tau)=
\end{equation}
$$=\sum_1^{\infty}\gl^{-n-1}\sum_0^{\infty}p_{\ell}(-\tl{\pp})\gl^{-\ell}
\pp_n\pp_1^2log\tau$$
and equate powers of $\gl$.  Write then ($f=\pp_1log\tau$)
\bq\label{33}
(\sum_1^{\infty}p_{km}\gl^{-k-m}f)(1+\sum_{n=1}^{\infty}\sum_{\ell=0}^{\infty}
p_{\ell}\pp_n\gl^{-n-\ell-1}f)=\sum\sum p_{\ell}\pp_n\gl^{-n-\ell-1}\pp_1f
\end{equation}
\bq\label{34}
\sum_{k,m=1}^{\infty}p_{km}\gl^{-k-m}f+\sum_{k,m=1}^{\infty}\sum_{n=1}^{\infty}
\sum_{\ell=0}^{\infty}p_{k,m}fp_{\ell}\pp_nf\gl^{-k-m-\ell-n-1}=
\end{equation}
$$=\sum_{n=1}^{\infty}\sum_{\ell=0}^{\infty} p_{\ell}\pp_n\gl^{-n-\ell-1}\pp_1f$$
\bq\label{35}
\gl^{-2}:\,\,p_{11}f=\pp_1^2f;\,\,\gl^{-3}:\,\,2p_{12}f=\pp_2\pp_1f+p_1\pp_1^2f;
\end{equation}
$$\gl^{-4}:\,\,(2p_{13}+p_{22})f+p_{11}f\pp_1f=\pp_3\pp_1f+p_1\pp_2\pp_1f+
p_2\pp_1^2f$$
This is ``doable" but becomes a bit tedious so we omit further terms.
\\[3mm]\indent
We can also encode matters in an elegant manner based on \cite{c4}.  Thus write L
in \eqref{31} as
\bq\label{36}
\sum_1^{\infty}p_kp_m\gl^{-k-m}log\tau=
\sum_1^{\infty}F_{km}\gl^{-k-m}=\sum_1^{\infty}(\sum_{k+m=j}F_{km})\gl^{-j}=
\sum_2^{\infty}Z_j\gl^{-j}
\end{equation}
(note $Z_1=0$).  Then following \cite{c4}, \eqref{31} becomes
\bq\label{37}
exp\sum_1^{\infty}Z_j\gl^{-j}=\sum_0^{\infty}p_i(Z_j)\gl^{-i}=1+
\sum_{n=1}^{\infty}\sum_{\ell=0}^{\infty}\gl^{-n-\ell-1}p_{\ell}(-\tl{\pp})
\pp_n\pp_1log\tau
\end{equation}
Note there is no $\gl^{-1}$ term on the right side which is balanced by
$Z_1=0$.
Then one arrives at
\bq\label{38}
\sum_2^{\infty}p_i(Z_j)\gl^{-i}=\sum_2^{\infty}\gl^{-i}\sum_{n+\ell=i-1}
p_{\ell}(-\tl{\pp})\pp_n\pp_1log\tau\Rightarrow
\end{equation}
$$\Rightarrow p_i(Z_j)=\sum_{n+\ell=i-1}p_{\ell}(-\tl{\pp})\pp_n\pp_1
log \tau$$
It is shown in \cite{t3} that the subset of Pl\"ucker relations involved in
\eqref{31} of
\eqref{36} is sufficient to determine the KP hierarchy.  Hence the Hirota
equations are encoded in \eqref{38} and can be expressed immediately in terms
of $u=\pp_1^2log\tau$ by writing
\bq\label{39}
Z_j=\sum_{k+m=j}\pp_1^2\tl{F}_{km};\,\,\tl{F}_{km}=\pp_1^2p_k(-\tl{\pp})
p_m(-\tl{\pp})log\tau=p_k(-\tl{\pp})p_m(-\tl{\pp})u
\end{equation}
Thus $F_{km}=\pp_1^2\tl{F}_{km}=F_{km}(u)$ and $Z_j=Z_j(u)$.  Similarly the
right side of the last equation in \eqref{38} can be written as $\pp^{-1}
\sum_{n+\ell=i-1}p_{\ell}(-\tl{\pp})\pp_nu$.  Consequently
\begin{theorem}
The classical KP Hirota equations can be written directly in terms of
$u=\pp^2log\tau$ via ${\bf (A26)}\,\,p_i(Z_j(u))=\sum_{n+\ell=i-1}\pp_1^{-1}
\pp_nu$.
\end{theorem}
\begin{corollary}
One can also compute using \eqref{32} as in \eqref{33}-\eqref{35}.  A standard
version of the Hirota equations in a new form can be given in terms of 
\eqref{26} and \eqref{27}.
\end{corollary}

\section{MORE ON QKP}
\renewcommand{\theequation}{4.\arabic{equation}}
\setcounter{equation}{0}

One reason for developing Hirota formulas in Section 3 was to be able to give
subsequently a q-version in order to write down the form of qKP equations
without explicitly compputing terms in an expansion ${\bf (A27)}\,\,
\pp_nL_q=[(L_q^n)_{+},L_q]$ or even $\pp_nQ=[Q^n_{+},Q]$ where $Q=D+a_0
+\sum_1^{\infty}a_nD^{-n}$ with $Df(x)=f(qx)$ (Frenkel form).  Formulas of this
type (e.g. originating in {\bf (A27)})
could be useful in order to compare with equations arising as in Examples 2.1,
2.2 in a q-context.  First recall \eqref{2}
$exp_q(xz)=exp(\sum_1^{\infty}c(xz)_i)= exp\sum (1-q)^k(xz)^k/k(1-q^k)$ and
(cf. \eqref{4})
$\psi_q=Se_q(xz)exp(\xi)=(1+\sum\tl{w}_iz^{-i})exp_q(xz)exp(\xi)$ with (cf.
$\eqref{15}) \psi_q=\tau_q(x,t-[z^{-1}])exp_q(xz)exp(\xi)/\tau_q$.  
One sort of expects $\psi_q(x,t)=\psi(t+[x]_q)$ since $exp_q(xz)exp(\xi)=exp
(\sum (t_iz^i+c(xz)_i)=exp(\sum(t_i+c(x)_i)z^i$ (recall $c(x)\sim [x]_q$) - see
below).  Thus \eqref{3} $\sum \tl{p}_k(x,t)z^k=e_q(xz)exp(\xi)$ implies
\bq\label{40}
\sum \tl{p}_k(x,t_i)z^k=exp\sum (t_i+c(x)_i)z^i=\sum
p_k(t_i+c(x)_i)z^k\Rightarrow \tl{p}_k(x,t)=p_k(t+c(x))
\end{equation}
Similarly $\psi^*_q(x,t)=\psi^*(t+[x]_q)$.  So what about formulas like \eqref
{226} or \eqref{26} or \eqref{27}?  We have ($\pp_1\sim \pp/\pp t_1$)
\bq\label{41}
\psi_q(x,t)=X_q(t,z)\tau_q(t)/\tau_q(t)=e_q(xz)e^{\xi}exp(-\sum_1^{\infty}
(\pp_i/i)z^{-i})\tau_q/\tau_q=\frac{e_q(xz)e^{\xi}\tau_q(t-[z^{-1})}{\tau_q(t)}
\end{equation}
since ${\bf
(A28)}\,\,exp(-\sum_1^{\infty}(\pp_i/i)z^{-i})\tau_q=\tau_q(t-[z^{-1}])=\tau
(t+c(x)-[z^{-1}])$. 
Now following \cite{i6} one looks at ${\bf (A29)}\,\,S=\sum_0^{\infty}
\tl{w}_i(t+[x]_q)D_q^{-i}$ with 
$(S^{-1})^*_{x/q}=\sum_0^{\infty}\tl{w}_j^*
(t+[x]_q)(-1)^jD^{-j}_{1/q}$).  Note from \eqref{41}
\bq\label{44}
\psi_q=S(e_q(xz)e^{\xi})=(1+\sum\tl{w}_jz^{-j})e_q(xz)e^{\xi}=e_q(xz)e^{\xi}\sum
p_{\ell}(-\tl{\pp})z^{-\ell}\tau_q/\tau_q\Rightarrow
\end{equation}
$$\Rightarrow 1+\sum\tl{w}_jz^{-j}=\frac{\sum
p_{\ell}(-\tl{\pp})z^{-\ell}\tau_q}{\tau_q}\Rightarrow \tl{w}_j=\frac
{p_j(-\tl{\pp})\tau_q}{\tau_q}\Rightarrow S=\sum_0^{\infty}
\frac{p_j(-\tl{\pp})\tau_q}{\tau_q}D_q^{-j}$$
as in \eqref{18} and ${\bf (A30)}\,\,\tl{w}_j=w_j(t+c(x))$.
As for Hirota recall {\bf (A9)}; however $t'\ne t$ is not specified here and
we will get more mileage from \eqref{7} in any event.
Thus we go to \eqref{7} for Hirota with $D_q$ terms.  Write $X_q\tau_q=e_q(xz)
exp(\xi)exp(-\sum\tl{\pp}_iz^{-i})\tau_q=e_q(xz)exp(\xi)\sum p_n(-\tl{\pp})
z^{-i}\tau_q$.  Also 
$\tl{X}_q=e_{1/q}(-xz)X(-t,z)$ so
$$\tl{X}_q\tau_q=e_{1/q}(-xz)exp(-\xi)\sum p_m(\tl{\pp})z^{-i}\tau_q$$  
Then
\eqref{7} becomes ($n>m$)
\bq\label{56}
0=\oint D^ne_q(xz)e^{\xi(t,z)}\sum
p_{\ell}(-\tl{\pp})z^{-\ell}\tau(t+c(x))\cdot
\end{equation}
$$D^{m+1}e_{1/q}(-xz)e^{-\xi(t',z)}\sum
p_k(\tl{\pp})z^{-k}\tau(t'+c(x))dz$$
For $n=1,\,m=0$ one has $De_q$ and $De_{1/q}$, with $D_qe_q=ze_q$ and 
for $e_{1/q}$ one goes to Remark 2.3.  Thus $D_{1/q}=qD_qD^{-1}$ so since
$D_{1/q}e_{1/q}(-xz)=-ze_{1/q}(-xz)$ there results
$qD_qe_{1/q}(-xz/q)=D_{1/q}e_{1/q}(-xz)=-ze_{1/q}(-xz)$.  Thus
$qD_qe_{1/q}(-xz)=D_{1/q}e_{1/q}(-xzq)=-qze_{1/q}(-xzq)=-qzDe_{1/q}(-xz)$
leading to $D_qe_{1/q}=-zDe_{1/q}$. We cannot directly replace D by $D_q$ in
\eqref{56} of course but for $n=m+1$ one can write $D^ne_q(xz)=e_q(xq^nz)$ and
$D^ne_{1/q}(-xz)=e_{1/q}(-xq^nz)=e_q(xq^nz)^{-1}$.  Hence for $n=m+1$
\eqref{56} becomes
\bq\label{57}
0=\oint e^{\xi(t,z)}\sum p_{\ell}(-\tl{\pp})z^{-\ell}D^n\tau_q
e^{-\xi(t',z)}\sum p_k(\tl{\pp})z^{-k}D^n\tau_q
\end{equation}
and this can be put in bilinear form as in \eqref{25}.  We state this as
\begin{proposition}
For $n=m+1$ \eqref{56} becomes \eqref{57} which can be put into bilinear form
as in \eqref{25}.  The resulting Hirota equations essentially ignore the $D^n$
term and produce standard formulas for $\tau(t+c(x))$ 
or $D^n\tau(t+c(x))$ in the $t$ variables.
\end{proposition}
\indent
Generally however we want to put $D_q$ explicitly into the act and for this we
need expressions for D in terms of $D_q$ (analogous to \eqref{122}).  One has
${\bf (A31)}\,\,D=1+(q-1)xD_q$ and we note
\bq\label{59}
qxD_qf=\frac{qxf(qx)-xf(x)}{(q-1)x}-f(x)\frac{(q-1)x}{(q-1)x}\Rightarrow
D_qx=qxD_q+1
\end{equation}
Consequently
\bq\label{60}
D^2=(1+(q-1)xD_q)(1+(q-1)xD_q)=1+2(q-1)xD_q+(q-1)^2xD_qxD_q=
\end{equation}
$$=1+(q^2-1)xD_q+q(q-1)^2x^2D_q^2$$
We note from \eqref{122}
$$\nm_q=\frac{(q^n-1)\cdots(q^{n-m+1}-1)}{(q^m-1)\cdots(q-1)}$$
\bq\label{61}
D^2=qx^2(q-1)^2D_q^2+(q^2-1)xD_q+1
\end{equation}
Similarly ${\bf (A32)}\,\,D_qx^2=q^2x^2D_q+[(q^2-1)/(q-1)]x$ and
\bq\label{62}
D^3=(q-1)^3q^3x^3D_q^3+qx^2(q-1)(q^3-1)D_q^2+(q^3-1)xD_q+1
\end{equation}
Another term should reveal the pattern but we prefer to use \eqref{122}
(cf. below).  In any event one will have a formula
\bq\label{63}
D^n=\sum_0^na_{nk}D_q^k
\end{equation}
with $a_{nk}$ to be determined.  One could go to 
\eqref{122} in the form (cf. $(\bgs)$)
\bq\label{69}
q^{n(n-1)/2}x^n(q-1)^nD_q^n=\sum_0^n(-1)^mq^{m(m-1)/2}\nm_qD^{n-m}
\end{equation}
but the calculations are unnecessarily complicated.
However a lovely formula was pointed out to the author by J. Cigler, based on his
article $(\bgs)$ Monatshefte f\"ur Mathematik, 88 (1979), 87-105, and we give the
result here with his proof.
\begin{theorem}
({\bf Cigler})
The $a_{nk}$ in \eqref{63} can be written in the form
\bq\label{702}
D^n=\sum_0^n{\nm}_q(q-1)^mx^mq^{m(m-1)/2}D_q^m
\end{equation}
\end{theorem}
To prove this one begins with the formula
\bq\label{701}
x^n=\sum_0^n{\nm}_q(x-1)(x-q)\cdots (x-q^{m-1})\Rightarrow 
\end{equation}
$$\Rightarrow D^n=\sum_0^n
{\nm}_q(D-1)(D-q)\cdots (D-q^{m-1})$$
(note here e.g. $x^m(x^{-1};q)_m=(x-1)(x-q)\cdots (x-q^{m-1})$).
Then to show \eqref{702}
one checks that $(D-1)(D-q)\cdots (D-q^{m-1})=(q-1)^m
q^{m(m-1)/2}x^mD_q^m$ and this is immediate via verification on $x^r$ for
$r=0,1,2,\cdots$.
\\[3mm]\indent
Now write in \eqref{56}
\bq\label{58}
0=\oint e_q(xq^nz)e^{\xi(t,z)}\sum p_{\ell}(-\tl{\pp})z^{-\ell}D^n\tau_q\cdot
\end{equation}
$$e_{1/q}(-xq^{m+1}z)e^{-\xi(t',z)}\sum
p_k(\tl{\pp})z^{-k}D^{m+1}\tau_qdz$$
Let $t'\to t+y$ and $t\to t-y$ so $e^{\xi}e^{-\xi'}\to exp(-2\sum y_iz^i)$ and
then we are faced with powers of z in $e_q(xq^nz)$ and $e_{1/q}(-xq^{m+1}z)=
1/e_q(xq^{m+1}z)$.  Here one could write for example
$e_q(xz)=exp(\sum_1^{\infty} [(1-q)^kx^kz^k/k(1-q^k)]$ and $n-m-1=s$ so that
\bq\label{588}
log D^ne_q(xz)-log D^{m+1}e_q(xz)=\sum_1^{\infty}\frac{(1-q)^kx^kz^k}{k(1-q^k)}
q^{(m+1)k}(q^{sk}-1)=\sum_1^{\infty}b_k^{ms}(q)x^kz^k
\end{equation}
This would introduce a term ($b_k=b_k^{ms}$)
\bq\label{589}
exp(\sum_1^{\infty}b_k^{ms}(q)x^kz^k=\sum_0^{\infty}p_{\nu}(b_k)(xz)^{\nu}
\end{equation}
into the integrand of \eqref{58}. 
For completeness we indicate the calculations involving \eqref{590}.  Thus
from \eqref{58} we
get
\bq\label{64}
\oint \sum_0^{\infty}p_{\nu}(b_k)x^{\nu}z^{\nu}\sum_0^{\infty}p_r(-2y)z^r
\sum_0^{\infty}p_{\ell}(-\tl{\pp})z^{-\ell}D^n\tau_q(t-y)\sum_0^{\infty}
p_k(\tl{\pp})z^{-k}D^{m+1}\tau_q(t+y)dz=0
\end{equation}
We cannot use the bilinear form now unless $n=m+1$ (cf. \eqref{57}) but
nevertheless one can write down a general residue formula. However the presence
of the
$\nu$ index complicates things horribly and produces an infinite number of
terms in each equation since all $\nu$ come into play.  
\\[3mm]\indent
{\bf REMARK 4.1}
It might be better 
to go directly to $D_q$ via
\eqref{63} and use the rules derived after \eqref{56}, namely (recall
$qDD_q=D_qD$ and $D_qD^{-1}=q^{-1}D_{1/q}$)
\bq\label{590}
D_q^ne_q(xz)=z^ne_q(xz);\,\,D_qe_{1/q}(-xz)=-zDe_{1/q}(-xz)\Rightarrow
\end{equation}
$$D_q^ne_{1/q}(-xz)=(-1)^nz^nq^nD^ne_{1/q}(-xz)$$
Note also
\bq\label{591}
D_q\left(\frac{1}{e_q}\right)=\frac{e_q-De_q}{e_qDe_q}=-\frac{D_qe_q}{e_qDe_q}
=-\frac{z}{De_q}
\end{equation}
The problem is to find a nice expression for $D^ne_q(xz)D^{m+1}e_{1/q}(-xz)
=D^ne_q/D^{m+1}e_q$ and \eqref{588} gives one version; some calculations in this
direction using interaction with $D_q$ are possible and we will discuss this
further in a forthcoming paper.  Another gambit is offered via
\cite{i6} (second paper) where the Frenkel approach for KP is recast in terms of
$\gD_qf(x)=f(qx)-f(x)=(D-1) f(x)$.  Thus
$$\gD_q^nf=\sum_0^{\infty}\nk(D^{n-k}\gD_q^kf)\gD_q^{n-k}$$
for $n\in {\bf Z}$.  Defining
\bq\label{ugh}
E_q(x,z)=exp\left(\frac{log(x)log(1+z)}{log(q)}\right)
\end{equation}
there results $\gD_qE_q=zE_q$ and $\gD_q^*E_{1/q}=zE_{1/q}$ where $(\sum
a_i\gD_q^i)^* =\sum \gD_{1/q}^ia_i$.  Then qKP takes the form
$L=\gD_q+a_0+\sum_1^{\infty} a_i\gD_q^{-i}$ with $\pp_iL=[L_{+}^j,L]$ and
$\psi_q=E_q(x,z)exp(\xi)
\tau_q(x,t-[z^{-1}])/\tau_q(x,t)$; here $\tau_q(x,t)=\tau(t+[[x]]_q)$ where
$[[x]]_q=((-1)^{k-1}log(x)/klog(q))$ for $k\geq 1$.  Also formulas relating
$D_q^n$ and $\gD_q^k$ should be simpler than those involving $D^n$ and $D_q^k$
given below in Theorem 4.2.
$\hfill\bs$
\\[3mm]\indent
\indent
{\bf REMARK 4.2.}
In keeping with Remark 2.1 we should determine now a formula for u in terms of 
$D_q$.  Note for qKdV the formula in Remark is $u=D_q\pp_1(log\tau_q(x,t)
\tau_q(qx,t))$ so $\pp_1$ is admitted.  Thus we leave $\pp_1$ and $t_1$ intact
even though it is tempting to abandon them.  To obtain a u formula we use
\eqref{44} where 
\bq\label{78}
S=1-\pp_1log\tau_qD_q^{-1}+\cdots
\end{equation} 
and ${\bf (A33)}\,\,L=SD_qS^{-1}$ or better $LS=D_qS$.  We recall first from
\eqref{11} 
\bq\label{79}
D_qf=DfD_q+(D_qf);\,\,D_q^{-1}f=D^{-1}fD_q^{-1}+D^{-2}D_qfD_q^{-1}+\cdots
\end{equation}
Then consider
\bq\label{80}
(D_q+a_0+\sum a_iD_q^{-i})(1+\tl{w}_1D_q^{-1}+\tl{w}_2D_q^{-2}+\cdots)=
\end{equation}
$$=(1+\tl{w}_1D_q^{-1}+\tl{w}_2D_q^{-1}+\cdots)D_q$$
There results
$$D_q+(D_q\tl{w}_1)D_q^{-1}+(D\tl{w}_1)+(D\tl{w}_2)D_Q^{-1}+(D_q\tl{w}_2)D_q^{-2}
+\cdots
+a_0+a_0\tl{w}_1D_q^{-1}+a_1D_q^{-1}+\cdots=$$
\bq\label{81}
=D_q+\tl{w}_1+\tl{w}_2D_q^{-1}+\cdots
\end{equation}
leading to
\bq\label{82}
\tl{w}_1=a_0+(D\tl{w}_1);\,\,(D_q\tl{w}_1)+(D\tl{w}_2+a_0\tl{w}_1+a_1=\tl{w}_2
\end{equation}
Consequently
\bq\label{83}
a_0=\tl{w}_1-(D\tl{w}_1);\,\,a_1=\tl{w}_2-(D\tl{w}_2)-\tl{w}_1^2+
\tl{w}_1D\tl{w}_1-D_q\tl{w}_1
\end{equation}
\begin{theorem}
Instead of $u=\pp^2log\tau$ as in KP one has a more complicated formula
for qKP, namely
\bq\label{84}
u=a_1=(1-D)\left(\frac{(1/2)(\pp_1^2-\pp_2)\tau_q}{\tau_q}\right)
-\left(\frac{\pp_1
\tau_q}{\tau_q}\right)^2+\frac{\pp_1\tau_q}{\tau_q}D\left(\frac{\pp_1\tau_q}
{\tau_q}\right)+D_q\frac{\pp_1\tau_q}{\tau_q}
\end{equation}
Note for $q\to 1$ this reduces to $\pp(\pp\tau/\tau)=\pp^2log\tau$ as desired.
The formula can be rewritten as
\bq\label{88}
u=-(q-1)xD_q\left(\frac{p_2(-\tl{\pp})\tau_q}{\tau_q}\right)+\frac{\pp_1\tau_q}
{\tau_q}(q-1)xD_q\left(\frac{\pp_1\tau_q}{\tau_q}\right)+D_q\left(
\frac{\pp_1\tau_q}{\tau_q}\right)
\end{equation}
\end{theorem}

\section{ZERO CURVATURE}
\renewcommand{\theequation}{5.\arabic{equation}}
\setcounter{equation}{0}

We go back now to Examples 2.1 and 2.2 and ask whether one can produce qKP
type equations from a first order differential calculus (FODC) of some sort.
One can (and we will later do this) of course look at differential calculi on
a quantum plane or q-deformed Heisenberg algebra which would produce $D_q$ or
$D_{q^2}$ operators in all variables.  First however let us try to link $D_q$
in the x variable with y and t as ``ordinary" variables subject to various 
noncommutativity relations among the differentials.  Consider first the
Burger's equation as in Example 2.1.
\begin{example}
Try e.g. $xdx=q^2dxx$ ($\gG_{+}$ style as in \cite{cxx,k1}) with $[dx,t]=adt,
[dx,dt]=0=(dt)^2=(dx)^2$, and $[dt,x]=adt$ for consistency.  Then $dx^n=
[[n]]_{q^{-2}}x^{n-1}dx=[[n]]_{q^2}dxx^{n-1}$; we write this here as
$\pp_qx^ndx$ (recall
$[[n]]_{q^2}=(q^{2n}-1)/(q^2-1)$).  Further take $dtt=tdt$ so 
\bq\label{186}
df=\pp_qfdx+(a+1)f_tdt
\end{equation}
Try $A=wdt+udx$ as before so
\bq\label{187}
dA=(\pp_qw+(a+1)u_t)dxdt
\end{equation}
For $A^2$ we note that $xt=tx$ is compatible with the conditions above and
e.g. ${\bf (A34)}\,\,wdtwdt=w(x,t)w(x+a)(dt)^2=0$ while for $udxudx$ one notes
$dxf(x)=D^{-2}fdx$ while ${\bf (A35)}\,\,dxt^m=t^mdx+mat^{m-1}dt$ implies that
$dxu=D_x^{-2}udx+aD_x^{-2}u_tdt$.  This leads to
\bq\label{189}
A^2=[w(x,t)w(x+a,t)+uaD_x^{-2}u_t+uD_x^{-2}w]dtdx
\end{equation}
Hence $dA+A^2=F=0$ requires
\bq\label{190}
w(x,t)w(x+a,t)+auD_x^{-2}u_t+uD_x^{-2}w+\pp_qw+(a-1)u_t=0
\end{equation}
Note if $q\to 1$ one obtains ${\bf (A36)}\,\,w^2+uw+w_x+u_t=0$ which for
$w=u_x$ implies ${\bf (A37)}\,\,u_x^2+uu_x+u_{xx}+u_t=0$
which is a kind of perturbation of the Burger's equation 
by a $(u_x)^2$ term (constants can be adjusted by changes of variables etc.)
$\hfill\bs$
\end{example}
\begin{example}
Consider next the q-plane situation
\bq\label{192}
(dx)^2=(dt)^2=0;\,\,xt-qtx=0;\,\,dxdt=-q^{-1}dtdx;\,\,xdx=q^2dxx;
\end{equation}
$$xdt=qdtx+(q^2-1)dxt;\,\,tdx=qdxt;\,\,tdt=q^2dtt$$
Then
\bq\label{193}
df=D_t^{-1}\pp_q^xfdx+\pp_q^tfdt
\end{equation}
(with $\pp_q$ as in Example 5.1).  Then $A=udx+wdt$ yields
\bq\label{194}
dA+A^2=0\leadsto -q\pp_q^tu+D_t^{-1}\pp_q^xw+uD_x^{-2}D_t^{-2}w-
\end{equation}
$$-\frac{wqt}{x}[D_x^{-1}D_t^{-1}w-D_x^{-3}D_t^{-1}w]-qwD_x^{-1}D_t^{-2}u=0$$
For $q\to q^{-1}$ evidently $\pp_qf=[f(q^{-2}x)-f(x)]/(q^{-2}-1)x\to
[f(q^2x)-f(x)]/(q^2-1)x$ and we write this latter as $\hat{\pp}_q$.  This gives
\bq\label{195}
-q^{-1}\hat{\pp}_q^tu+D_t\hat{\pp}_q^xw+uD_x^2D_t^2w-q^{-1}wD_xD_t^2u-
\frac{wt}{qx}[D_xD_tw-D_x^3D_tw]=0
\end{equation}
Note that as $q\to 1$
\eqref{195} becomes
\bq\label{196}
-\pp_tu+\pp_xw+w-wu=0
\end{equation}
Taking $w=u_x$ one has 
\bq\label{197}
-u_t+u_{xx}+u_x-uu_x=0\leadsto u_t+uu_x+u_{xx}-u_x=0
\end{equation}
which is again a kind of perturbed Burger's equation with perturbation $-u_x$.
$\hfill\bs$
\end{example}
\indent
{\bf REMARK 5.1.}
What should a q-Burger's equation look like
and how do we produce it by zero curvature conditions as indicated?  Based on
Example 2.1 one might consider a candidate equation as
\bq\label{198}
\hat{\pp}_q^tw+a\hat{\pp}_q^xw^2+b(\hat{\pp}_q^x)^2w=0
\end{equation}
with suitable factors of $D_x,\,\, D_t$ inserted.  
Consider then in \eqref{195} a substitution $w=q\hat{\pp}_q^xu$; this leads to
\bq\label{102}
-q^{-1}\hat{\pp}_q^tu+aD_t(\hat{\pp}_q^x)^2u-aq^{-1}\hat{\pp}_q^xu(D_xD_t^2u)=
-aD_x^2D_t^2\hat{\pp}_q^xu+
\end{equation}
$$+\frac{ta^2}{qx}(\hat{\pp}_q^xu)D_xD_t(1-D_x^2)\hat{\pp}_q^xu$$
which is a perturbation of \eqref{198} for example involving powers of
$\hat{\pp}_q^xu$.$\hfill\bs$
\\[3mm]\indent
{\bf REMARK 5.2.}
We note also the Cole-Hopf transformation for Burger's equation.  Thus given
${\bf (A38)}\,\,u_t+2uu_x+\ga u_{xx}=0$ one can write
\bq\label{103}
u=-\ga\frac{\psi_x}{\psi}=-\ga\pp(log(\psi));\,\,u_x=-\ga\pp^2(log(\psi));\,\,
u_{xx}=-\ga\pp^3(log(\psi))
\end{equation}
Then ${\bf (A39)}\,\,-\ga\pp_t\pp_xlog(\psi)+\ga^2(\psi_x/\psi)^2-\ga^2\pp^3
log(\psi)=0$.  But 
\bq\label{104}
\pp\left(\frac{\psi_x}{\psi}\right)-\pp^3log(\psi)=-\pp\left(\frac{\psi_{xx}}
{\psi}\right)
\end{equation}
Consequently for $\ga<0$
\bq\label{105}
\ga\pp_t\pp_xlog(\psi)-\ga^2\pp(\psi_{xx}/\psi)=0\Rightarrow\psi_t-\ga\psi_{xx}
=c\psi (=0\,\,for\,\,c=0)
\end{equation}
If there is an additional term $\gb u_x$ in {\bf (A38)} then the same method
gives rise to ${\bf (A40)}\,\,\psi_t-\ga\psi_{xx}=\gb\psi_x$ which can be
treated with the same facility as the case $\gb=0$.  In any event one reduces
matters to heat equations and a q-Cole-Hopf transformation should
exist.$\hfill\bs$
\\[3mm]\indent
{\bf REMARK 5.3.}
These examples may yet go somewhere but perhaps suffer from the imposed
partially canonical structure based on $xdx=q^2dxx$ (cf. \cite{cxx,k1}).
Let us be more simple minded now and work from the qkDv model where
\bq\label{200}
L^2=D^2_q+(q-1)xuD_q+u;\,\,u=D_q\pp_1log[\tau(x,t)D\tau(x,t)]
\end{equation}
We will broach this matter later using zero curvature ideas for KdV.  For now
consider situations related to Burger's equation and try to see what are some
good models for q-differential equations.  For KdV (cf. \cite{t2}) the zero
curvature condition is ${\bf (A41)}\,\,\pp_mB_n-\pp_nB_m-[B_n,B_m]=0$ 
based on $\pp_nL^2=[B_n,L^2]$ where 
$B_n=L^n_{+}$.  Here $L=\pp_q+s_0 +s_1\pp_q^{-1}+\cdots$ and ${\bf (A42)}\,\,
u_1=(q-1)xu=s_0+D(s_0)$.  One knows ${\bf
(A43)}\,\,\pp_tu=(\pp_qu)-(\pp_q^2s_0) -(\pp_qs_0^2)$ where $\pp_q\sim D_q$
here.  Further from {\bf (A42)} one has
$(q-1)xu=(q-1)x\pp_1\pp_q(log\tau+log(D\tau))=s_0+Ds_0$ which in fact implies
${\bf (A44)}\,\,s_0=(q-1)x\pp_1\pp_qlog\tau$ (note $qDD_q=D_qD$ so $D(x\pp_1
\pp_qlog\tau)=x\pp_1qDD_qlog\tau=x\pp_1\pp_qDlog\tau$).  Hence in principle we
can compute the qKdV equation via {\bf (A43)} but it is complicated.  Thus 
(recall $\pp_q(fg)=Df\pp_qg+(\pp_qf)g$)
\bq\label{201}
\pp_qs_0=(q-1)\pp_1[qx\pp_q^2log\tau+\pp_qlog\tau];\,\,\pp_q^2s_0=
\end{equation}
$$=(q-1)\pp_1[q^2x\pp_q^3log\tau+q\pp_q^2log\tau+\pp_q^2log\tau]$$
and
$$s_0^2=(q-1)^2x^2(\pp_1\pp_qlog\tau)^2;\,\,\pp_qs_0^2=(q-1)^2[q^2x^2
\pp_q(\pp_q\pp_qlog\tau)^2+
\frac{q^2-1}{q-1}x(\pp_1\pp_qlog\tau)^2];$$
\bq\label{202}
\pp_q(\pp_1\pp_qlog\tau)^2=[D(\pp_1\pp_qlog\tau)]\pp_q(\pp_q\pp_qlog\tau)+
\pp_q(\pp_1\pp_qlog\tau)(\pp_1\pp_qlog\tau)
\end{equation}
This makes it appear feasible only to use the bilinear or some other Hirota form
to deal with qKdV.  The expression {\bf (A43)} is not good in terms of u.  We
will return to these matters later.$\hfill\bs$

\end{document}